\pdfoutput=1
\RequirePackage{ifpdf}
\ifpdf 
\documentclass[pdftex]{sigma}
\else
\documentclass{sigma}
\fi

\numberwithin{equation}{section}

\numberwithin{theorem}{section}
\numberwithin{proposition}{section}
\numberwithin{lemma}{section}
\numberwithin{corollary}{section}
\numberwithin{definition}{section}
\numberwithin{example}{section}
\numberwithin{remark}{section}
\numberwithin{note}{section}
\numberwithin{conjecture}{section}

\newtheorem{question}{Question}
\numberwithin{question}{section}

\def \QQ {\mathcal{Q}}

\def \U {\mathcal{U}}

\def \C {\mathbb{C}}
\def \F {\mathbb{F}}
\def \K {\mathbb{K}}
\def \P {\mathbb{P}}
\def \Q {\mathbb{Q}}
\def \R {\mathbb{R}}

\def \Z {\mathbb{Z}}

\def \Gm  {\mathrm{G}_m}

\def \geq {\geqslant}

\def \leq {\leqslant}

\newcommand{\ang}[1] {\langle #1 \rangle}

\def \kappa {\varkappa}

\def\={\;=\;}
\def\bal{\begin{aligned}}
\def\eal{\end{aligned}}

\newcommand{\Spec}{{\text{Spec }}}

\DeclareMathOperator{\Hom}{Hom}

\DeclareMathOperator{\Aut}{Aut}

\DeclareMathOperator{\rk}{rank}

\begin{document}

\allowdisplaybreaks

\renewcommand{\thefootnote}{$\star$}

\renewcommand{\PaperNumber}{005}

\FirstPageHeading

\ShortArticleName{Upper Bounds for Mutations of Potentials}

\ArticleName{Upper Bounds for Mutations of Potentials\footnote{This
paper is a contribution to the Special Issue ``Mirror Symmetry and Related Topics''. The full collection is available at \href{http://www.emis.de/journals/SIGMA/mirror_symmetry.html}{http://www.emis.de/journals/SIGMA/mirror\_{}symmetry.html}}}

\Author{John Alexander CRUZ MORALES~$^{\dag^1}$ and Sergey GALKIN~$^{\dag^2\dag^3\dag^4\dag^5}$}

\AuthorNameForHeading{J.A.~Cruz Morales and S.~Galkin}

\Address{$^{\dag^1}$ Department of Mathematics and Information Sciences, Tokyo Metropolitan University, \\
\hphantom{$^{\dag^1}$}~Minami-Ohsawa 1-1, Hachioji, Tokyo 192-037, Japan}
\EmailDD{\href{mailto:cruzmorales-johnalexander@ed.tmu.ac.jp}{cruzmorales-johnalexander@ed.tmu.ac.jp}, \href{mailto:alekosandro@gmail.com}{alekosandro@gmail.com}}

\Address{$^{\dag^2}$~Kavli Institute for the Physics and Mathematics of the Universe, The University of Tokyo,\\
\hphantom{$^{\dag^2}$}~5-1-5 Kashiwanoha, Kashiwa, 277-8583, Japan}

\Address{$^{\dag^3}$~Independent University of Moscow, 11 Bolshoy Vlasyevskiy per., 119002, Moscow, Russia}

\Address{$^{\dag^4}$~Moscow Institute of Physics and Technology, 9 Institutskii per.,\\
\hphantom{$^{\dag^4}$}~Dolgoprudny, 141700, Moscow Region, Russia}
\EmailDD{\href{mailto:Sergey.Galkin@phystech.edu}{Sergey.Galkin@phystech.edu}}

\Address{$^{\dag^5}$ Universit\"at Wien, Fakult\"at f\"ur Mathematik, Garnisongasse 3/14, A-1090 Wien, Austria}

\ArticleDates{Received May 31, 2012, in f\/inal form January 16, 2013; Published online January 19, 2013}

\Abstract{In this note we provide a new, algebraic proof of the excessive Laurent phenomenon for mutations of potentials (in the sense of [Galkin~S., Usnich~A., {P}reprint IPMU 10-0100, 2010]) by introducing to this theory the analogue of the upper bounds from~[Berenstein~A., Fomin~S., Zelevinsky~A., \textit{Duke Math.~J.} \textbf{126} (2005), 1--52].}


\Classification{13F60; 14J33; 53D37}

\renewcommand{\thefootnote}{\arabic{footnote}}
\setcounter{footnote}{0}

\section{Introduction}

The idea of mutations of potentials was introduced in \cite{GU}
and the Laurent phenomenon was established in the two dimensional case
by means of birational geometry of surfaces.
More precisely, in op.\ cit.\ the authors considered
a toric surface $X$ with a rational function $W$ (a~potential),
and using certain special birational transformations (mutations),
they established the (excessive) Laurent phenomenon which roughly says that
if $W$ is a Laurent polynomial whose mutations are Laurent polynomials, then
all subsequent mutations of these polynomials are also Laurent polynomials
(see Theorem~\ref{ulemma} in Appendix~\ref{appendix2} for a precise statement of the excessive Laurent phenomenon
as established in~\cite{GU}).
The motivating examples of such potentials
come from the mirror images of special Lagrangian tori
on del Pezzo surfaces~\cite{FOOO}
and Auroux's wall-crossing formula relating invariants of dif\/ferent tori~\cite{AU}.

The cluster algebras theory of Fomin and Zelevinsky~\cite{FZ} provides
an inductive way to construct \emph{some} birational transformations
of $n$ variables as a consecutive composition of elementary ones (called elementary mutations)
with a choice of $N=n$ directions at each step.

The theory developed in \cite{GU} can be seen as an extension of the theory of cluster algebras~\cite{FZ} when the number of directions of mutations $N$ is allowed to be (much) bigger than the number of variables $n$, but at least one function remains to be a Laurent polynomial after all mutations. So, it is natural to try to extend the machinery of the theory of cluster algebras for this new setup. The main goal of this paper is to give the f\/irst step in such an extension by means of the introduction of the upper bounds (in the sense of~\cite{BFZ}) and establishing the excessive Laurent phenomenon \cite{GU} in terms of them.
It is worth noticing that a further generalization can be done and in a forthcoming work~\cite{CG} we plan to study the quantization of the mutations of potentials and their upper bounds. Naturally, this quantization can be seen as an extension of the theory of quantum cluster algebras developed in~\cite{BZ,BZ1} and the theory of cluster ensembles in~\cite{FG}.

The upper bounds introduced in this paper can be described as a collection of regular functions that remain regular after one elementary mutation in any direction. Thus, we can establish the main result of this paper in the following terms (see Theorem~\ref{ublemma} for the exact formulation).

\begin{theorem*}[Laurent phenomenon in terms of the upper bounds]
The upper bounds are preserved by mutations.
\end{theorem*}

Aside from providing a new proof for the excessive Laurent phenomenon and the already mentioned generalization in the quantized setup, the algebraic approach that we are introducing here is helpful for tackling the following two problems:
\begin{enumerate}\itemsep=0pt
\item
Develop a higher dimensional theory (i.e.\ dimension higher than $2$) for the mutations of potentials. Some work in that direction is carried out in \cite{CGA}.
\item
Present an explicit construction to compactify Landau--Ginzburg models (Problem~44 of~\cite{GU}).
\end{enumerate}

In the present paper we do not deal with the above two problems (only a small comment on~2 will be made at the end of the paper). We plan to give a detailed discussion of them in~\cite{CG} too. We just want to mention that the new algebraic approach has interesting geometrical applications.

Some words about the organization of the text are in order. In Section~\ref{sec-def} we extend the theory developed in~\cite{GU} to lattices of arbitrary rank and general bilinear forms (i.e., we can consider even degenerate and not unimodular forms) and introduce the notion of upper bounds in order to establish our main theorem. In Section~\ref{section3} we actually establish the main theorem and present its proof when the rank of the lattice is two and the form is non-degenerate which is the case of interest for the geometrical setup of~\cite{GU}. In the last section some questions and future developments are proposed. For the sake of completeness of the presentation we include two appendices. In Appendix~\ref{sec-bfz} we review some def\/initions of~\cite{BFZ} and brief\/ly compare their theory with ours. Appendix~\ref{appendix2} is dedicated to presenting the Laurent phenomenon in terms of~\cite{GU}.

\section{Mutations of potentials and upper bounds} \label{sec-def}

Now we present an extension of the theory of mutations of potentials~\cite{GU}
(as formulated by the second author and Alexandr Usnich)
and introduce our modif\/ied def\/initions with the new def\/inition of upper bound.
Notice that a slightly dif\/ferent theory (which f\/its into the framework of this paper, but not \cite{GU})
is used in our software code\footnote{\url{http://member.ipmu.jp/sergey.galkin/degmir.gp}.}.

\subsection{Combinatorial data}

Let $(\cdot,\cdot) : L^* \times L \to \Z$
be the canonical pairing between
a pair of dual lattices
$L \simeq \Z^r$ and $L^* = \Hom(L,\Z) \simeq \Z^r$.

In what follows the lattice $L$ is endowed with a skew-symmetric bilinear integral form $\omega: L \times L \to \Z$
(we use the notation $\ang{v,v'} = \omega(v,v')$).
In the most important (both technically, and from the point of view of applications) case
$r = \rk L = 2$, we have $\Lambda^2 L \simeq \Z$,
so all integer skew-symmetric bilinear forms are integer multiples $\omega_k = k \omega_1$ ($k \in \Z$)
where a~generator~$\omega_1$ is f\/ixed by the choice of orientation on $L \otimes \R$
so that $\omega_1((1,0),(0,1)) = 1$.
We would occasionally use notations
$\ang{\cdot,\cdot}_1 = \omega_1(\cdot,\cdot)$
and
$\ang{\cdot,\cdot}_k = \omega_k(\cdot,\cdot)$.

The bilinear form $\omega$ gives rise to a map $i = i_{\omega} : L \to L^*$ that sends an element $v \in L$ into a linear form
$i_\omega(v) \in L^*$ such that $(i_\omega(v),v') = \omega(v,v')$ for any $v' \in L$.
The map $i_\omega$ is an isomorphism $\iff$ the form $\omega$ is non-degenerate and unimodular,
when~$\omega$ is non-degenerate but not unimodular the map $i$ identif\/ies the lattice $L$
with a full sublattice in $L^*$ of index~$\det \omega$,
f\/inally if $\omega$ is degenerate then both the kernel and the cokernel of the map $i_\omega$ has positive rank.

We would like to have some functoriality,
so we consider a category whose objects are given by pairs $(L,\omega)$ of the lattice $L$
and a skew-symmetric bilinear form $\omega$,
and the morphisms $\Hom((L',\omega'),(L,\omega))$ are linear maps $f : L' \to L$
such that $\omega' = f^* \omega$,
i.e.\ $\omega(v_1,v_2) = \omega'(f(v_1),f(v_2))$ for all $v_1,v_2 \in L'$.
Any linear map $f : L' \to L$ def\/ines an adjoint $f^* : L^* \to L'^*$
and if it respects the bilinear forms,
then $i_{\omega'} = f^* i_{\omega} f$.

For a vector $u \in L$ we def\/ine
a \emph{symplectic reflection} $R_u$
and a \emph{piecewise linear mutation} $\mu_u$
to be the (piecewise)linear automorphisms of the set $L$ given by the formulae
\begin{gather*}
R_{\omega,u} (v) = v + \omega(u,v) u, \\
\mu_{\omega,u} v = v + \max(0,\omega(u,v)) u.
\end{gather*}
For any morphism $f \in \Hom((L',\omega'),(L,\omega))$ and any vector $u \in L'$ we have
$R_{\omega, fu} f = f R_{\omega',u}$ and
$\mu_{\omega, fu} f = f \mu_{\omega',u}$.
Indeed, $f \mu_u v = f (v + \max(0,\omega'(u,v)) u) = fv + \max(0,\omega'(u,v)) (fu) =  fv + \max(0,\omega(fu,fv)) (fu) = \mu_{fu} (fv)$.

Note that
$R_{a \omega, bu} = R_{\omega,u}^{a b^2}$ for all $a,b \in \Z$
and
$\mu_{a \omega, bu} = \mu_{\omega,u}^{a b^2}$ for all $a,b \in \Z_+$.
However
$\mu_{\omega,-u} v = - \mu_u (-v) =  v + \min(0,\omega(u,v)) u$,
hence
$\mu_{\omega,-u} \mu_{\omega,u} = R_{\omega,u}$.
Both $R_{\omega,u}$ and $\mu_{\omega,u}$ are invertible:
$R_{\omega,u}^{-1} v = R_{-\omega,u} v  = v - \omega(u,v) u$,
$\mu_{\omega,u}^{-1} v = \mu_{-\omega,-u} v = R_{\omega,u}^{-1} \mu_{\omega, -u} v = v - \max(0,\omega(u,v)) u$.
Note that
$\mu_{\omega,-u}^{-1}  v = R_{\omega,-u}^{-1} \mu_{\omega,u}  v = R_{\omega,u}^{-1} \mu_{\omega,u}(v) = v - \min(0,\omega(u,v)) u$.
Therefore, changing $\max$ by $\min$ and $+$ by $-$, simultaneously,
corresponds to changing the form $\omega$ to the opposite $-\omega$.
Further we omit $\omega$ from the notations of $R_u$ and $\mu_u$ where the choice of the form is clear.

The underlying combinatorial gadget of our story
is a collection of $n$ vectors in $L$:
\begin{definition}
An \emph{exchange collection} $V$ is an element of $L^n$, i.e.\ an $n$-tuple $(v_1,\dots,v_n)$ of
vectors~$v_i \in L$.
Some~$v_i$ may coincide.
For a~vector $v$ its \emph{multiplicity} $m_V(v)$ in the exchange collection $V$ equals the number of vectors in $V$ that coincide with $v$:
$m_V (v) = \# \{ 1 \leq i \leq n : v_i = v\}$. We say that an exchange collection $V'$ is a \emph{subcollection} of exchange collection $V$ if
$m_{V'} \leq m_V$.
Equivalently, one may def\/ine an exchange collection $V$ by its (non-negative integer) multiplicity function
$m_V : L \to \Z_{\geq 0}$. In this case $n = \sum\limits_{v \in L} m_V(v)$.
\end{definition}

The exchange collections could be pushed forward by morphisms $f \in \Hom((L',\omega'),(L,\omega))$:
$v_1',\dots,v_n' \in L'^n$ will go to $f v_1,\dots, f v_n \in L^n$.
This gives rise to a natural diagonal action of $\Aut (L,\omega) = \operatorname{Sp}(L,\omega)$
on $L^n$. This action commutes with the permuting action of~$S_n$.

A vector $n \in L$ is called \emph{primitive}
if it is nonzero and its coordinates are coprime, i.e.~$n$ does not belong to the sublattice $k L$ for any $k>1$,
in other words $n$ is not a multiple of other vector in~$L$.
We denote the set of all primitive vectors in~$L$ as~$L_1$.
Similarly one can def\/ine primitive vectors in the dual lattice~$L^*$.
Note that if $\det \omega \neq \pm{1}$
then $i_\omega(n)$ may be a non-primitive element of~$L^*$ even for primitive elements~$n \in L_1$.

\subsection{Birational transformations}
Consider the group ring $\Z[L^*]$~-- ring of Laurent polynomials of $r$ variables.
Its spectrum $T = \Spec \Z[L^*] \simeq \Gm^r(\Z)$ is the $r$-dimensional torus over the integers,
in particular $T(\C) = \Hom(L^*,\C^*)$,
$L^* = \Hom(T,\Gm)$ is the lattice of characters of $T$
and
$L = \Hom(\Gm,T)$ is the lattice of $1$-parameter subgroups in~$T$.
Def\/ine \emph{the ambient field} $\K = \K_L = \QQ(L^*)$ as the fraction f\/ield of $\Z[L^*]$
extended by all roots of unity ($\QQ = \Q(\exp(2 \pi i \Q)$).

A vector $u \in L$ def\/ines a birational transformation of $\K_L$ (and its various subf\/ields and subrings) as follows
\begin{equation*}
\mu_{u,\omega} : \ X^m \to X^m \big(1+X^{i_\omega(u)}\big)^{(u,m)}.
\end{equation*}

If $f : T_1 \to T_2$ is a rational map between two tori, and $u : \Gm \to T_1$ is a one-parameter subgroup of $T$
then its image $f u : \Gm \to T_2$ is not necessarily a one-parameter subgroup, but asymptotically behaves like one,
this def\/ines a \emph{tropicalization} map $T(F) : \Hom(\Gm,T_1) \to \Hom(\Gm, T_2)$. The tropicalization of the birational
map $\mu_{u,\omega} : T_1 \to T_2$ is the piecewise-linear map
$\mu_{u,\omega} : L_1 \to L_2$ def\/ined in the previous subsection.

One can easily see most of the relations of the previous subsection on the birational level.
For example, $\mu_{-u} \mu_u = \mu_u \mu_{-u} = R_u$ and $R_v \mu_u R_v^{-1} = \mu_{R_v u}$,
where $R_u$ is the homomorphism of the torus $T$ given by
$R_{u,\omega} : X^m \to X^{m+(u,m) i_\omega u}$. Also $R_{a u, b \omega} = R_{u,\omega}^{a^2 b}$ for any $a,b \in \Z$,
and $\mu_{a u,\omega} = (\mu_{u, a \omega})^a$ for any $a \in \Z$, however neither of them is a power of $\mu_{u,\omega}$.\footnote{Since $(1+x^a)$ is not a power of $(1+x)$.} In particular, $(\mu_{u,\omega})^{-1} = \mu_{-u,-\omega}$.

Note that if $M \subset L^*$ is some sublattice of $L^*$ that contains $i_\omega(u)$ then $\mu_u$ preserves the fraction
f\/ield of $\Z[M] \subset \Z[L^*]$.
For any morphism $f \in \Hom((L',\omega'), (L,\omega))$ and a vector $u \in L'$ we have a homomorphism
$f^* : \Z[L^*] \to \Z[L'^*]$ and two birational transformations
$\mu_u \in \Aut \K_{L'}, \mu_{fu} \in \Aut \K_L$ that commute:
$\mu_u f^* = f^* \mu_{fu}$.

\begin{remark} \label{fun1}
We have the following functoriality of the mutations with respect to the lattice~$L$:
let $L' \subset L$ be a sublattice of index $k$ in the lattice~$L$,
so $L^* = \Hom(L,\Z)$ is a sublattice of index~$k$ in~$L'^* = \Hom(L',\Z)$,
and assume that the vector $u$
lies in the sublattice $L'$.
Then the Abelian group $G = (L/L')$ of order $k$ acts on $\QQ[L'^*]$,\footnote{An element $n$ in $L$ multiplies monomial $X^{m'}$ by the root of unity $\exp((2 \pi i) (n,m'))$, here $(n,m')$ is bilinear pairing between
$L$ and $L'^*$ with values in $\Q$ extended by linearity from the pairing $L' \otimes L'^* \to \Z$.}
and its invariants is the subring $\QQ[L^*]$,
so~$G$ acts on the torus $T' = \Spec \QQ[L'^*]$ and the torus $T = \Spec \QQ[L^*]$ is the quotient-torus $T = T' / G$,
let $\pi : T' \to T$ be the projection to the quotient.
The vector $u$ def\/ines the birational transformation~$\mu_{u,T}$ of the torus $T$ and the birational transformation~$\mu_{u,T'}$ of the torus~$T'$.
Then the mutation~$\mu_u$ commutes with the action of the group~$G$ and with the projections:
$\pi \mu_{u,T'} = \mu_{u,T} \pi$
and $g \mu_{u,T'} = \mu_{u,T'} g$
for any $g \in G$.
\end{remark}

\subsubsection{Rank two case}
Let us see the mutations explicitly in case $\rk L = 2$.
Let $e_1$, $e_2$ be a base of $L$ and $f_1$, $f_2$ be the dual base of~$L^*$,
so $(e_i,f_j) = \delta_{i,j}$. Also let $x_i = X^{f_i}$ be the respective monomials in $\Z[L^*]$.
For
the skew-symmetric bilinear form $\omega_k$ def\/ined by $\omega_k(e_1,e_2) = k$
and a vector $u = u_1 e_1 + u_2 e_2 \in L$
we have
$i_{\omega_k} (u_1 e_1 + u_2 e_2) = (-k u_2) f_1 + (k u_1) f_2$ and so
\begin{equation*}
\mu_{u,\omega_k} : \ (x_1,x_2) \to \big(x_1 \cdot \big(1+x_1^{-k u_2} x_2^{k u_1}\big)^{u_1}, x_2 \cdot \big(1+x_1^{-k u_2} x_2^{k u_1}\big)^{u_2}\big),
\end{equation*}
in particular the inverse map to $\mu_{u,\omega_1}$ is given by
$\mu_{-u,-\omega_1} : (x_1,x_2) \to (x_1 \cdot (1+x_1^{u_2} x_2^{-u_1})^{-u_1}, x_2 \cdot (1+x_1^{u_2} x_2^{-u_1})^{-u_2})$.
In particular, $\mu_{(0,1)}^* f = f\big(x_1,\frac{x_2}{1+x_1}\big)$.

For any matrix $A = \left(\begin{matrix} a & b  \\ c & d\end{matrix} \right) \in \operatorname{SL}(2,\Z) = \operatorname{Sp}(L,\omega)$
there is a regular automorphism of the torus $t_A^* (x_1,x_2) = (x_1^a x_2^b, x_1^c x_2^d)$.
Conjugation by this automorphism acts on the set of mutations:
$\mu_{A u}^* = (t_A \mu_u t_A^{-1})^*$.
So any mutation commutes with an inf\/inite cyclic group given by the stabilizer of $u$ in $\operatorname{Sp}(L,\omega)$,
explicitly if $u=(0,1)$ then in coordinates $(x_1,x_2)$ and $(x_1,x_2' = x_1 x_2)$ the mutation $\mu_{(0,1)}$ is given by the same formula.
Also every mutation commutes with $1$-dimensional subtorus of~$T$,
in case of $u=(0,1)$ the action of the subtorus is given by $(x_1,x_2) \to (x_1, \alpha x_2)$.

\subsection{Mutations of exchange collections and seeds}
Let $L$ be a lattice equipped with a bilinear skew-symmetric form $\omega$.
A~\emph{cluster} $\mathbf{y} \in \K_L^m$ is a collection
$\mathbf{y} = (y_1,\dots,y_m)$
of $m$ rational functions
$y_i \in \K_L$.
We call $\mathbf{y}$ \emph{a base cluster} if $\mathbf{y} = (y_1,,\dots,y_r)$
is a base of the ambient f\/ield $\K_L$.
A~$C$-seed (supported on  $(L,\omega)$)
is a pair $(\mathbf{y},V)$ of a cluster $\mathbf{y} \in \K_L^m$ and an exchange collection
$V = (v_1,\dots,v_n) \in L^n$.
A~$V$-seed (supported on  $(L,\omega)$)
is a pair $(W,V)$ of a rational function $W \in \K_L$ and
an exchange collection $V = (v_1,\dots,v_n) \in L^n$.

Given two exchange collections
 $V'=(v_1',\dots,v_n') \in L'^n$
and
 $V = (v_1,\dots,v_n) \in L^n$
we say that $V'$ is \emph{a mutation} of $V$
in the direction $1\leq j\leq n$ and denote it by $V' = \mu_j V$ if
under the given identif\/ication $s_j : L \simeq L'$
we have $v'_j = s_j(-v_j)$
and
$v'_i = s_j(\mu_{v_j} v_i)$ for $i\neq k$.

The mutation of a $C$-seed $(\mathbf{y},V)$ in the direction $1\leq j\leq n$
is a new $C$-seed $(\mathbf{y}_j,V_j)$ where
$V_j = \mu_j V$ is a mutation of the exchange collection,
and $\mathbf{y}_j = \mu_{v_j,\omega} \mathbf{y}$ where each variable is
transformed by the birational transformation $\mu_{v_j,\omega}$.

The identity $\mu_{-u} \mu_u = R_u$ implies that $\mu_j(\mu_j(V))$ and $V$ are
related by the $\operatorname{Sp}(L,\omega)$-trans\-for\-ma\-tion~$R_u$.

\subsection[Upper bounds and property $(V)$]{Upper bounds and property $\boldsymbol{(V)}$}

\begin{definition}[property $(V)$] \label{pv}
We say a $V$-seed $(W,V)$ satisf\/ies property $(V)$ if~$W$ is a Laurent polynomial and for all $v \in L$ the functions ${(\mu_{v}^*)}^{m_V(v)} W$ are also Laurent polynomials.
\end{definition}

In this paper we introduce the upper bound of an exchange collection.
\begin{definition}[upper bounds] \label{ub}
For a $C$-seed $\Sigma = (\mathbf{y},V)$
def\/ine its \emph{upper bound} $\U(\Sigma)$ to be the $\QQ$-subalgebra of $\K_L$ given by
\begin{equation*}
\U(\Sigma) = \QQ\big[\mathbf{y}^{\pm 1}\big] \cap \big(\cap_{v \in L} \QQ\big[{(\mu_{v}^*)}^{m_V(v)} \mathbf{y}^{\pm 1}\big]\big).
\end{equation*}
In case $\mathbf{y}$ is a base cluster (by abuse of notation)
we denote $\U(\Sigma)$ just by~$\U(V)$.
\end{definition}

The upper bounds def\/ined here are a straightforward generalization of the upper bounds in~\cite{BFZ},
but also they can be thought of as the gatherings of all potentials satisfying property~$(V)$.

\begin{proposition}[relation between property $(V)$ and upper bounds] \label{vub1}
The upper bound $\U(V)$ of an exchange collection $V$ consists of all functions $W \in \K_L$
such that the $V$-cluster $(W,V)$ satisfies property~$(V)$.
\end{proposition}

\begin{proposition} \label{prop-isom}
Any morphism $f : (L,\omega) \to (L',\omega')$
induces a dual morphism $f^* : L'^* \to L^*$,
a homomorphism of algebras $f^* : \Z[L'^*] \to \Z[L^*]$.
Assume that this homomorphism has no kernel\footnote{One can bypass this assumption by def\/ining the upper bound $\U(L,\omega;V)$ as a subalgebra in
some localization of $\Z[L^*]$ determined by the exchange collection $V$.}.
Then it induces a homomorphism of upper bounds $f^* : \U(L',\omega';f V) \to \U(L,\omega;V)$.
In particular, if $f$ is an isomorphism, then maps $f^*$ and $(f^{-1})^*$ establish the isomorphisms
between the upper bounds $f^* : \U(L',\omega';f V) \simeq \U(L,\omega;V)$.
\end{proposition}

\begin{proposition} \label{fun2}
Consider a seed $\Sigma = (L,\omega;v_1,\dots,v_n)$.
For a sublattice $L' \subset L$ that contains all vectors $v_i \in L' \subset L$
consider the seed $\Sigma' = (L',\omega_{|_{L'}};v_1,\dots,v_n)$.
By Remark~{\rm \ref{fun1}} there is a~natural action of
$G = L/L'$ on $\K_{L'}$
with $\K_L = \K_{L'}^G$.
Moreover,
the action of $G$ obviously preserves
the property of being a~Laurent polynomial
$($i.e.\ it preserves the subalgebras $\QQ[L'^*])$,
and the mutations~$\mu_{v_i}$ commute with the $G$-action.
Thus the upper bound with respect to the overlattice $L$
is the subring of $G$-invariants of the upper bound with respect to the sublattice~$L'$:
$\U(\Sigma) = \U(\Sigma')^G = \U(\Sigma') \cap \QQ[L^*]$.
\end{proposition}

\section{Laurent phenomenon}\label{section3}

In what follows we restrict ourselves to the case $\rk L = 2$, $\omega$ is a non-degenerate form and the vectors of exchange collection are primitive,
however none of these conditions is essential.

Next theorem is the analogue of Theorem~1.5 in~\cite{BFZ}, presented here as Theorem~\ref{bfz-theorem}.

\begin{theorem}[Laurent phenomenon in terms of upper bounds] \label{ublemma}
Consider two $C$-seeds: $\Sigma = (L,\omega;v_1,\dots,v_n)$ and
$\Sigma' = (L',\omega';v_1',\dots,v_n')$.
If $\Sigma' = \mu_i \Sigma$ is a mutation of $\Sigma$ in direction $1\leq i\leq n$
then
the upper bounds for $\Sigma$ and $\Sigma'$ coincide: $\U(\Sigma) = \mu_{v_i}^* \U(\Sigma')$.
As a~corollary, if a seed~$\Sigma'$ is obtained from a seed $\Sigma$
by a sequence of mutations, then the upper bound $\U(\Sigma')$
equals to the upper bound~$\U(\Sigma)$ under identification of the ambient field
by composition of the birational mutations.
\end{theorem}

By Proposition \ref{vub1} Theorem~\ref{ublemma} is equivalent to the next corollary,
which is easier to check in practice
and has almost the same consequences as the main theorem
of~\cite{GU}, presented here as Theorem~\ref{ulemma}.
\begin{corollary}[$V$-lemma] \label{vlemma}
If $V$-seeds $\Sigma$ and $\Sigma'$ are related by a mutation
then the seed $\Sigma$ satisfies property $(V)$ $\iff$
the seed $\Sigma'$ satisfies property~$(V)$.
\end{corollary}

In the rest of this section we prove Theorem~\ref{ublemma}.
Our proof is quite similar to that of~\cite{BFZ}\footnote{See Appendix~\ref{sec-bfz} and Remark~\ref{compare-bfz} for the detailed comparison.}:
The set-theoretic argument
reduces the problem to
exchange collection $V$ with small number of vectors~($1$~or~$2$) without counting of multiplicities.
Actually, when the collection~$V$ has only one vector
the equality of the upper bounds is obvious from the def\/initions.
When the exchange collection consists of two base vectors one can explicitly compute the upper bounds
and compare them. Finally, the case of two non-base non-collinear vectors is thanks to functoriality.

First of all, let us f\/ix the notations.
If the rank two lattice~$L$ is generated by a pair of vectors~$e_1$ and~$e_2$,
then the dual lattice $L^* = \Hom(L,\Z)$ has the dual base $f_1$, $f_2$ determined by
$(f_i,e_j) = \delta_{i,j}$. The form $\omega$ is uniquely determined by its value
$k = \omega(e_1,e_2)$, and further we denote this isomorphism class of forms by~$\omega_k$.
We assume that $k\neq 0$, i.e.\ the form $\omega$ is
non-degenerate\footnote{If $k=0$ then $\omega=0$ and all mutations are trivial.},
by swapping~$e_1$ and~$e_2$ one can exchange $k$ to $-k$.
A base~$e_i$ of $L$ corresponds to a base $x_i = X^{f_i}$ of $\Z[L^*]$.

\begin{lemma} \label{l41a}
Let $V$ be an exchange collection in $(L,\omega)$ and $\Sigma = (L,\omega;V)$ be the respective seed.
\begin{enumerate}\itemsep=0pt
\item[$1.$]
If $V$ is empty, then obviously $\U(L,\omega;V) = \QQ[L^*]$.
\item[$2.$]
Otherwise, let $V_\alpha$ be a set of exchange collections
such that for any $v\in L$ we have
$m_V(v) = \max_{\alpha} m_{V_\alpha}(v)$.
Then
\[ \U(V) = \cap_{\alpha} \U(V_\alpha) . \]
\item[$3.$]
In particular, if for a vector $v \in L$ we define $V_v = m_V(v) \times v$
to be an exchange collection that consists of a single vector $v$
with multiplicity $m_V(v)$ and $\Sigma_v = (L,\omega;V_v)$ be the respective seed,
then
$ \U(\Sigma) = \cap_{v \in L} (\QQ[\mathbf{y}^\pm] \cap \QQ[{(\mu_v^*)}^{m_V(v)} \mathbf{y}^\pm])
= \cap_{v \in L} \U(\Sigma_v)$.
In other words, the upper bound of a $C$-seed $\Sigma = (L,\omega;\mathbf{y},V)$ can be expressed as the intersection of
the upper bounds for its $1$-vector subseeds.
\item[$4.$]
Let $V$ consist of a vector $v_1$ with multiplicity $m_+ \geq 1$,
a vector $v_2 = - v_1$ with multiplicity \mbox{$m_- \geq 0$},
and vectors $v_k$ $(k\geq 3)$ that are non-collinear to $v_1$ with some multiplicities \mbox{$m_k \geq 0$}.
Consider exchange subcollections $V_0 = \{m_+\times v_1, m_-\times (-v_1)\}$
and $V_k = \{1\times v_1,$ $m_k \times v_k\}$ $(k\geq 3)$.
Then
\[ \U(V) = \U(V_0) \cap \U(V_3) \cap \U(V_4) \cap \cdots . \]
\item[$5.$]
Let $V' = \mu_1 V$ be an exchange collection obtained by mutation of~$V$ in~$v_1$;
it consists of vector $-v_1$ with multiplicity $m_- + 1\geq 1$,
vector $v_1$ with multiplicity $m_+ - 1 \geq 0$
and vectors $v_k' = \mu_{v_1} v_k$ $(k\geq 3)$ with multiplicities $m_k$.
Similarly to the previous step define $V_0' = \{(m_- + 1)\times (-v_1), (m_+ - 1)\times v_1\}$
and $V_k' = \{1\times (-v_1), m_k \times v_k'\}$ $(k \geq 3)$. Then
\[ \U(V') = \U(V_0') \cap \U(V_3') \cap \U(V_4') \cap \cdots. \]
\item[$6.$]
Hence, to proof Theorem~{\rm \ref{ublemma}} it is necessary and sufficient
to show that
\[ \U(V_0) = \mu_{v_1}^* \U(V_0')
\qquad \text{and} \qquad
\U(V_k) = \mu_{v_1}^* \U(V_k') \quad(\text{for all} \ k\geq 3 ). \]
We will prove these equalities in Proposition~{\rm \ref{l43c1}} and Lemma~{\rm \ref{l46}}.
\end{enumerate}
\end{lemma}

\begin{proposition} \label{isom2d}
Let $v_1,v_2 \in L$ be a pair of vectors $v_1 = a e_1 + b e_2$, $v_2 = c e_1 + d e_2$
such that  $ad-bc = 1$.
Consider the lattice $L'$ with the base $e_1'$, $e_2'$
and the form $\omega'(e_1',e_2') = \omega(v_1,v_2)$;
let~$f_1'$,~$f_2'$ be the dual base of~$L'^*$.
Consider a map $m : L \to L'$ given by
$m(e_1) = d e_1' - b e_2', m(e_2) = -c e_1' + a e_2'$;
note that
$m(v_1) = m(a e_1 + b e_2) = e_1'$ and  $m(v_2) = m(c e_1 + d e_2) = e_2'$.
The dual isomorphism $m^*: L'^* \to L^*$ is given by the transposed map
$m^* (f_1') =  d f_1 - c f_2$ and $m^*(f_2') = -b f_1 + a f_2$.
Let $z_1 = X^{f_1'} = x_1^d x_2^{-c}$ and $z_2 = X^{f_2'} = x_1^{-b} x_2^a$.
Since map $m^*$ is invertible by Proposition~{\rm \ref{prop-isom}} it gives the equality
\[ \U(L,\omega;m_1\times v_1,m_2\times v_2) = \U\big(L',\omega';m_1\times e_1',m_2\times e_2'\big)\big|_{{z_1 = x_1^d x_2^{-c},\,z_2 = x_1^{-b} x_2^a}}. \]
\end{proposition}

\begin{lemma} \label{l42}
Assume a seed $\Sigma = (L,\omega;m_1\times v_1)$ consists of a unique vector $v_1$ with multiplicity $m_1 \geq 1$.
\begin{enumerate}\itemsep=0pt
\item[$1.$]
If $v_1 = e_2 = (0,1)$
then the upper bound $\U(\Sigma)$ consists of all Laurent polynomials $W$ of the form
$W = \sum_l c_l(x_1) x_2^l$
where $c_l \in \QQ[x_1^\pm]$ and for $l \leq 0$ we have that $c_l$ is divisible by
$(1+x_1^k)^{-m_1 l}$.
Moreover, $\U(\Sigma) = \QQ\big[x_1^\pm, x_2, \frac{1}{x_2^{\prime \prime}} = \frac{(1+x_1^k)^{m_1}}{x_2}\big]$.
\item[$2.$]
If $v_1 = a e_1 + b e_2 = (a,b)$ is an arbitrary primitive vector
then
$\U(\Sigma) = \QQ\big[z^{\pm}, z_1, \frac{(1 + z^k)^{m_1}}{z_1}\big]$
 where
$z =  \frac{x_1^a}{x_2^b}$,
$z_1 = x_1^r x_2^s$
and $(r,s)\in\Z^2$ satisfies $rb+sa=1$.
\end{enumerate}
\end{lemma}

\begin{proof}
Recall that mutation in the direction $e_2$ is given by
$x_1' = x$ and $x_2' = \frac{x_2}{1+x_1^k}$.
Assume we have a Laurent polynomial $W = \sum_{l \in \Z} c_l(x_1) x_2^l$.
Then $W$ can be expressed in terms of $x_1$ and $x_2'$ as
$W = \sum_l c_l(x_1) (1+x_1^k)^l (x_2')^l$.
This function is a Laurent polynomial in terms of $(x_1,x_2')$ $\iff$
$c_l(x_1) (1+x_1^k)^l$ is a Laurent polynomial of $x_1$ for all $l$. This is equivalent to
$c_l$ being divisible by $(1+x_1^k)^{-l}$ for $l \leq 0$.
Similarly if we do $m_1$ mutations then $x_2^{\prime \prime} = \frac{x_2}{(1+x_1^k)^a}$
and $W = \sum c_l (1+x_1^k)^{m_1 l} (x_2^{\prime \prime})^l$
so for $l\leq 0$ we have that $c_l$ is divisible by $(1+x_1^k)^{-m_1 l}$.
Let $c_l(x_1) = (1+x_1^k)^{-m_1l} c'_{-l}(x_1)$ for $l < 0$, $c'_l$ are also Laurent polynomials.
Denote $W_+ = \sum_{l \geq 0} c_l(x_1) x_2^l$ and $W_- = \sum_{l <0} c_l(x_1) x_2^l = \sum_{l > 0} c'_l(x_1) (x_2^{\prime \prime})^{-l}$.
Then obviously both $W_+$ and $W_-$ belong to $\QQ\big[x_1^\pm,x_2,\frac{1}{x_2^{\prime \prime}}\big]$.
The reverse inclusion is straightforward.

Part (2) follows from Proposition~\ref{isom2d}.
\end{proof}

\begin{proposition} \label{lhz2}
Let exchange collection $V$ consists of a vector $v_1 = e_2 = (0,1)$ with multiplicity $m_1\geq 0$
and its inverse $v_2 = -v_1 = -e_2 = (0,-1)$ with multiplicity $m_2 \geq 0$,
\begin{enumerate}\itemsep=0pt
\item[$1.$]
The upper bound $\U(L,\omega;m_1\times e_2,m_2\times (-e_2))$
consists of all Laurent polynomials $W$ of the form
$W = \sum_l c_l(x_1) x_2^l$ where $c_l \in \QQ[x_1^\pm]$ and for $l \leq 0$ we have that $c_l$ is divisible by $(1+x_1^k)^{-m_1l}$
and for $l \geq 0$ we have that $c_l$ is divisible by $(1+x_1^k)^{m_2l}$.
\item[$2.$]
$\U(L,\omega;m_1\times e_2,m_2\times(-e_2)) = \QQ\big[x_1^\pm, x_2(1+x_1^k)^{m_2},\frac{(1+x_1^k)^{m_1}}{x_2}\big].
$
\end{enumerate}
\end{proposition}

\begin{proof}
The f\/irst statement is a straightforward corollary of Lemmas~\ref{l41a} and~\ref{l42}(1).
The proof of the second statement is similar to the end of the proof of Lemma~\ref{l42}(2):
separate the Laurent polynomial $W$ into positive and negative parts~$W_+$ and~$W_-$;
then both parts lie in the ring $\QQ\big[x_1^\pm,x_2 (1+x_1^k)^{m_2},\frac{(1+x_1^k)^{m_1}}{x_2}\big]$.
\end{proof}

\begin{proposition} \label{l43c1}
Assume a seed $\Sigma$ consists of a vector $v_1 = (0,1)$ with multiplicity $m_1$
and its inverse $-v_1 = (0,-1)$ with multiplicity $m_2$.
Then its mutation $\Sigma' = \mu_1(\Sigma)$ consists of~$v_1$ and~$-v_1$ with respective multiplicities $m_1-1$ and $m_2+1$.
Then $\U(\Sigma) = \U(\Sigma')$.
\end{proposition}
\begin{proof}
By Proposition~\ref{lhz2} the upper bounds are expressed as:
$\U(\Sigma) = \QQ\big[x_1^\pm,x_2 (1+x_1^k)^{m_2}$, $\frac{(1+x_1^k)^{m_1}}{x_2}\big]$,
$\U(\Sigma') = \QQ\big[x_1'^\pm,\frac{(1+x_1'^k)^{m_1-1}}{x_2'},x_2' (1+x_1'^k)^{m_2+1}\big]$.
Since $x_1' = x_1$ and $x_2' = \frac{x_2}{1+x_1^k}$ we have the desired equality of the upper bounds.
\end{proof}

\begin{proposition} \label{lhz3}
Assume that the seed $\Sigma = (L,\omega;m_1\times v_1,m_2\times v_2)$
consists of vectors $v_1$ with multiplicity $m_1 \geq 0$
and $v_2$ with multiplicity $m_2 \geq 0$.
\begin{enumerate}\itemsep=0pt
\item[$1.$]
If $v_1 = e_1$ and $v_2 = e_2$,
then the upper bound $\U(\Sigma)$ equals
$\QQ\big[x_1, x_2, \frac{(1 + x_2^k)^{m_1}}{x_1}, \frac{(1 + x_1^k)^{m_2}}{x_2}\big]$.
\item[$2.$]
If $v_1 = a e_1 + b e_2$ and $v_2 = c e_1 + d e_2$ with $ad-bc = 1$
then the upper bound $\U(\Sigma)$ equals
$\QQ\big[z_1, z_2, \frac{(1 + z_2^k)^{m_1}}{z_1}, \frac{(1 + z_1^k)^{m_2}}{z_2}\big]$
with
$z_1 = x_1^dx_2^{-c}$
and
$z_2 = x_1^{-b}x_2^a$.
\end{enumerate}
\end{proposition}

\begin{proof}
For the f\/irst case,
by Lemmas \ref{l41a} and \ref{l42}
we have
$\U(\Sigma) = \QQ\big[x_1^\pm,x_2,\frac{(1+x_1^k)^{m_2}}{x_2}\big] \cap \QQ\big[x_2^\pm,x_1,\frac{(1+x_2^k)^{m_1}}{x_1}\big]$. If $m_1 = m_2 = 1$,
by Proposition~4.3 of~\cite{BFZ} (with $|b_{12}| = |b_{21}| = b = c = k$ and $q_1=q_2=r_1=r_2=1$)
this intersection equals
$\QQ\big[x_1,x_2,\frac{1+x_2^k}{x_1},\frac{1+x_1^k}{x_2}\big]$.
Lemma~\ref{l42} covers cases with $m_1 = 0$ or $m_2 = 0$.
If $m_1$ and $m_2$ are greater than $1$,
the proof of Proposition~4.3 in~\cite{BFZ} can be easily modif\/ied to include the case we need since $x_2 \frac{1 + x_1^k}{x_2}(1 + x_1^k)^{m_{2}-1} = (1+ x_1^k)^{m_2}$ and $x_1 \frac{1 + x_2^k}{x_1}(1 + x_2^k)^{m_{1}-1} = (1+ x_2^k)^{m_1}$ , then the intersection equals $\QQ\big[x_1,x_2,\frac{(1+x_2^k)^{m_1}}{x_1},\frac{(1+x_1^k)^{m_2}}{x_2}\big]$.

Part (2) follows from Proposition~\ref{isom2d}.
\end{proof}

\begin{lemma} \label{l46}
Let $\Sigma = (L,\omega;1 \times v_1, m_2 \times v_2)$ be a seed of two non-collinear vectors $v_1$ and $v_2$
with $m(v_1) = 1$ and $m(v_2) = m_2 \geq 0$
and
$\Sigma' = \Sigma_1 = (L'=L,\omega'=\omega; v_1' = -v_1, m_2 \times (v_2' = \mu_{v_1} v_2) )$
be the mutation of the seed $\Sigma$ in~$v_1$.
Then $\U(\Sigma) = \mu_{v_1}^* \U(\Sigma')$.\footnote{Denote $\Sigma_2 = \{(\mu_{v_2})^{m_2} v_1, -v_2 \}$~-- $m_2$-multiple mutation of $\Sigma$ in $v_2$,
and $\Sigma'_2 = \{ (\mu_{v_2'})^{m_2} v_1',  -v_2' \}$~-- $m_2$-multiple mutation of $\Sigma'$ in $v_2'$.
We are going to prove that
 $\U(\Sigma) = \QQ[\mathbf{y}] \cap \QQ[\mathbf{y}' = \mathbf{y_1}] \cap \QQ[\mathbf{y_2}]$
equals to
$\U(\Sigma') = \QQ[\mathbf{y}] \cap \QQ[{\mathbf{y}'} = \mathbf{y_1}] \cap \QQ[\mathbf{y_2'}]$.}
\end{lemma}
\begin{proof}
First of all note that, $\omega'(v_1',v_2') = -\omega(v_1,v_2)$ and since
$\mu_{-v_1} \mu_{v_1} = R_{v_1}$ it is suf\/f\/icient to consider only the case
$\omega(v_1,v_2) > 0$.
We f\/irst consider the case when $v_1 = e_1 = (1,0)$ and $v_2 = e_2 = (0,1)$;
denote $k = \omega(e_1,e_2) > 0$.
Let $e_1'$, $e_2'$ be the base of $L'$ that corresponds to~$e_1$,~$e_2$ under the natural
identif\/ication of $L' \simeq L$;
f\/inally consider a base $e_1''$, $e_2''$ of $L'$
given by $e_1'' = v_2' = k e_1' + e_2'$, $e_2'' = v_1' = -e_1'$.
Let $f_1'$, $f_2'$ and $f_1''$, $f_2''$ be the respective dual bases of~$L'^*$.
Thus we have one natural regular system of coordinates
$x_1 = X^{f_1}$, $x_2 = X^{f_2}$
on the torus $T = \Spec \Z[L^*]$,
and two regular systems of coordinates
$x_1' = X^{f_1'}$, $x_2' = X^{f_2'}$; $x_1'' = X^{f_1''}$, $x_2'' = X^{f_2''}$
on the torus $T' = \Spec \Z[L'^*]$.
Taking $a=k$, $b=1$, $c=-1$ and $d=0$
in Proposition~\ref{isom2d}
we have that $x_1'' = x_2'$ and $x_2'' = \frac{x_2'^k}{x_1'}$, since $x_2' = x_2$ and $x_1' = \frac{x_1x_2^k}{1+x_2^k}$ (they are the mutations of~$x_1$ and~$x_2$ with respect to $v_1$), thus what we need to show is that the rings
$\QQ\big[x_1,x_2,\frac{1+x_2^k}{x_1}, \frac{1+x_1^k}{x_2}\big]$
and $\QQ\big[x_1,x_2, \frac{1+x_2^k}{x_1}, \frac{x_1^{k} + (1 + x_2^k)^k}{x_1^kx_2}\big]$ are equal.
We will f\/irst show that
$\frac{x_1^{k} + (1 + x_2^k)^k}{x_1^kx_2} \in \QQ\big[x_1,x_2, \frac{1+x_2^k}{x_1}, \frac{1+x_1^k}{x_2}\big]$.
We have that $ \frac{x_1^{k} + (1 + x_2^k)^k}{x_1^kx_2} = \big(\frac{1 + x_1^k}{x_2}\big)\big(\frac{(1+x_2^k)^k}{x_1^k}\big) - \sum\limits_{j=1}^k \frac{k!}{j!(k-j)!} x_2^{kj-1}$. Clearly the expression in the right side belongs to $\QQ\big[x_1,x_2, \frac{1+x_2^k}{x_1}, \frac{1+x_1^k}{x_2}\big]$.
Now, we will show that $\frac{1+x_1^k}{x_2} \in \QQ\big[x_1,x_2, \frac{1+x_2^k}{x_1}, \frac{x_1^k + (1 + x_2^k)^k}{x_1^kx_2}\big]$. We have that $\frac{1+x_1^k}{x_2} = x_1^k  \frac{x_1^{k} + (1 + x_2^k)^k}{x_1^kx_2} - \sum\limits_{j=1}^k \frac{k!}{j(k-j)!} x_2^{kj-1}$. Again, clearly the expression in the right side belongs to $\QQ\big[x_1,x_2, \frac{1+x_2^k}{x_1}, \frac{x_1^{k} + (1 + x_2^k)^k}{x_1^kx_2}\big]$. Thus, we have the equality between the rings. Similarly, if the multiplicity of $v_2$ is $m_2 > 1$, we have that $\QQ\big[x_1,x_2,\frac{(1+x_2^k)}{x_1}, \frac{(1+x_1^k)^{m_2}}{x_2}\big] = \QQ\big[x_1,x_2, \frac{1+x_2^k}{x_1}, \frac{(x_1^{k} + (1 + x_2^k)^k)^{m_2}}{x_1^{m_2k}x_2}\big]$.\footnote{The argument for showing the equality of these two rings is the same of that when $m_2 = 1$, but the computations are slightly longer, so we omit them.}
If $v_1$ and $v_2$ are another basis of $\Z^2$ the result follows from
Proposition~\ref{isom2d}.
In case~$v_1$ and~$v_2$ is a pair of non-collinear vectors which are not a basis for~$\Z^2$,
consider the sublattice $L' \subset L$ generated by $e_1' = v_1$ and $e_2' = v_2$
with the form $\omega' = \omega|_{{L'}}$.
As we just saw upper bounds with respect to the sublattice coincide:
$\U(L',\omega';v_1, m_2\times v_2) = \mu_{v_1}^* \U(L',\omega';-v_1,m_2\times \mu_{v_1} v_2)$.
Now the statement follows from the Proposition~\ref{fun2}.
\end{proof}

\begin{remark}
If a mutation of a Laurent polynomial with integer coef\/f\/icients
happened to be a Laurent polynomial, then its coef\/f\/icients are also integer.
Let $u \in L$ be a primitive vector
and $W, W' \in \QQ[L^*]$ be a pair of Laurent polynomials with arbitrary coef\/f\/icients
such that $W = \mu_u^* W'$.
Then $W \in \Z[L^*] \iff W' \in \Z[L^*]$.
\end{remark}

\begin{proof}
Choose coordinates on $L$ so that $u = e_2$.
Assume $W$ has integer coef\/f\/icients.
By Lemma~\ref{l42}(1) $W = \sum c_l(x_1) x_2^l$
and $W' = \sum_{l\in\Z} c_l'(x_1) x_2'^l$
with $c_l' = c_l (1+x_1^k)^{-l}$ for all $l\in \Z$.
Clearly $W \in \Z[L^*] \iff$ all coef\/f\/icients of $W$ are integer $\iff$
for all $l\in \Z$ all coef\/f\/icients of $c_l$ are integer.
Since $(1+x_1^k) \in \Q[x_1^{\pm}]$ it is clear that $c_l' \in \Q[x_1^{\pm}]$.
Recall that for a Laurent polynomial $P \in \Q[x]$ its Gauss's content
$C(P) \in \Q$ is def\/ined as the greatest common divisor of all its coef\/f\/icients:
if $P = \sum a_i x^i$ then $C = \gcd(a_i)$.
Clearly $C(P) \in \Z \iff P \in \Z[x^{\pm}]$.
Gauss's lemma says that $C(P \cdot P') = C(P) \cdot C(P')$.
Since $C(1+x_1^k) = \gcd(1,1) = 1$
we see that $C(c_l') = C(c_l) \cdot 1^{-l} = C(c_l)$,
hence $c_l' \in \Z[x_1^\pm] \iff c_l \in \Z[x_1^\pm]$.
\end{proof}

\section{Questions and future developments}

In the introduction was pointed out that our def\/inition of upper bounds makes plausible to consider a quantum version of mutations of potentials and the corresponding quantum Laurent phenomenon. On the other hand, in \cite{K} a non-commutative version of the Laurent phenomenon is discussed. Thus, we would like to ask:

\begin{question}
Is it possible to consider a non-commutative version of the Laurent phenomenon for mutation of potentials and develop a theory of upper bounds in this context?
\end{question}

In \cite{GU} the following problem (Problem 44) was proposed

\begin{question} \label{ques2}
Construct a fiberwise-compact canonical mirror of a Fano variety as a~gluing of open charts given by $($all$)$ different toric degenerations.
\end{question}

Conjecture \ref{conj} (which will be proved in \cite{CG}) gives a partial answer for the above question.

\begin{conjecture} \label{conj}
For the $10$ potentials $W$ $($i.e., $(W_1,W_2,\dots,W_9,W_Q))$ listed in {\rm \cite{GU}} $($or rather the exchange collections $V$ $($resp.\ $V_1,\dots,V_9,V_Q))$ the upper bound $\U(V)$ is the algebra of polynomials in one variable. Moreover, this variable is $W$.
\end{conjecture}

Conjecture~\ref{conj} is useful for symplectic geometry as long as one knows two (non-trivial) properties of the FOOO's potentials~$m_0$~\cite{FOOO} (here $W = m_0$):
\begin{enumerate}\itemsep=0pt
\item
$W$ is a Laurent polynomial (this is some kind of convergence/f\/initeness property).

\item
$W$ is transformed according to Auroux's wall-crossing formula~\cite{AU},
and more specif\/ically by the mutations described in Section~\ref{sec-def}. The directions of the mutations/walls are encoded by an exchange collection~$V$.
\end{enumerate}

What we believe is that once one knows these assumptions, one should be able to prove that some disc-counting potential equals some particularly written $W$ (formally) without any actual disc counting. Needless to say this is a speculative idea.

\appendix
\section{Review of the classical cluster algebras, upper bounds\\ and Laurent phenomenon}
\label{sec-bfz}

In this appendix we review some results of the f\/irst section of~\cite{BFZ}:
approach to Laurent phenomenon via upper bounds by Berenstein, Fomin and Zelevinsky, and make a brief comparison between their theory and the one presented here. We will denote the framework of cluster algebras developed by Berenstein, Fomin and Zelevinsky in~\cite{BFZ} by BFZ.

\subsection{Def\/initions of exchange matrix, coef\/f\/icients, cluster and seed}

Fix $n$-dimensional lattice $L \simeq \Z^n$.
The underlying combinatorial gadget in the theory of cluster algebras is a $n\times n$ matrix.

\begin{definition}[exchange matrix~$B$]
An exchange matrix is a sign-skew-symmetric $n \times n$ integer matrix $B =(b_{ij})$:
for any $i$ and $j$, either $b_{ij} = b_{ji} = 0$ or $b_{ij} b_{ji} < 0$.
\end{definition}
Obviously a skew-symmetric matrix is sign-skew-symmetric,
and for simplicity we assume further that $B$ is skew-symmetric.

Any matrix $B$ can be considered as an element of $L^* \otimes L^*$.
Skew-symmetric matrices are then identif\/ied with $\wedge^2(L)$.

Let $\P$ be the \emph{coefficient group}~-- an Abelian group without torsion written multiplicatively.
Fix an ambient f\/ield $\F$ of rational functions on $n$ independent variables with coef\/f\/icients in (the f\/ield of fractions of)
the integer group ring $\Z\P$.

\begin{definition}[coef\/f\/icients]
A \emph{coefficient tuple} $\mathbf{p}$ is an $n$-tuple of pairs $(p_i^+,p_i^-) \in \P^2$.
\end{definition}

Finally the non-combinatorial object of the theory is a cluster.
\begin{definition}[BFZ-cluster]
A \emph{cluster} $\mathbf{x} = (x_1,\dots,x_n)$ is a transcendence basis of $\F$ over the f\/ield of fractions of $\Z\P$.
Let $\Z\P[\mathbf{x}^{\pm 1}]$ denote the ring of Laurent polynomials of $x_1,\dots,x_n$ with coef\/f\/icients in $\Z\P$.
\end{definition}

\begin{definition}[BFZ-seed]
A \emph{seed} (or BFS-seed) is a triple $(\mathbf{x},\mathbf{p},B)$
of a cluster, coef\/f\/icients tuple and exchange matrix.
\end{definition}

\begin{remark}[action of the symmetric group~$S_n$]
As noticed in~\cite{BFZ} the symmetric group~$S_n$ naturally acts
on exchange matrices, coef\/f\/icients, clusters,
and hence seeds
by permutating indices~$i$.
\end{remark}

\subsection{Mutations}

For each $1\leq k \leq n$ we can def\/ine the mutation of exchange matrix $B$,
of a pair $(B,\mathbf{p})$ and of a~seed $(B,\mathbf{p},\mathbf{x})$.

\begin{definition}[mutation $\mu_i$ of an exchange matrix~$B$] \label{mub}
Given an exchange matrix $B = (b_{ij})$ and an index $1 \leq k \leq n$
def\/ine $\mu_k B = B' = (b'_{ij})$ as follows:
$b_{ik}' = - b_{ik}$,
$b_{kj}' = - b_{kj}$,
and otherwise
$b_{ij}' = b_{ij} + \frac{|b_{ik}|b_{kj} + b_{ik}|b_{kj}|}{2}$.
\end{definition}
It is easy to check that $\mu_k(\mu_k(B)) = B$.

\begin{definition}[mutations of coef\/f\/icients] \label{mup}
Given an exchange matrix $B$ and a coef\/f\/icients tuple $\mathbf{p}$
def\/ine a mutation of the coef\/f\/icients in direction $k$ as any new $n$-tuple
$(p_i^{\prime+},p_i^{\prime-})$ that satisf\/ies
$\frac{p_i^{\prime+}}{p_i^{\prime-}} = (p_k^+)^{b_{ki}} \frac{p_i^+}{p_i^-}$ if $b_{ki} \geq 0$
and
$\frac{p_i^{\prime+}}{p_i^{\prime-}} = (p_k^-)^{b_{ki}} \frac{p_i^+}{p_i^-}$ if $b_{ki} \leq 0$.
\end{definition}

In this def\/inition the choice of a new $n$-tuple has $(n-1)$ degrees of freedom.
This ambiguity is not important, however one of the ways of curing this ambiguity
is by considering tuples with $p^-=1$. Also one can get rid of coef\/f\/icients by considering the trivial tuples $p^+=p^-=1$.

\begin{definition}[mutations of seeds]
The mutation of a seed $\Sigma = (\mathbf{x},\mathbf{p},B)$
in the direction $1 \leq k \leq n$ is a new seed $\Sigma' = (\mathbf{x}_k,\mathbf{p}',B')$
where $B' = \mu_k B$ is a mutation of the exchange matrix~$B$ in the direction~$k$,
$\mathbf{p}'$ is a mutation of $\mathbf{p}$ using $B$ in the direction $k$ (Def\/inition~\ref{mup}),
and~$\mathbf{x'}$ is def\/ined as follows:
$x_k' x_k = P_k(\mathbf{x}) = p_j^+ \prod\limits_{b_{ik}>0} x_i^{b_{ik}} + p_k^- \prod\limits_{b_{ik}<0} x_i^{-b_{ik}}$
and $x'_i = x_i$ for $i\neq k$.
\end{definition}

The next def\/inition is a technicality required by \cite{BFZ} for the proof.
\begin{definition}
A seed $\Sigma$ is called \emph{coprime} if polynomials $P_1,\dots,P_n$ are pairwise coprime in~$\Z\P[\mathbf{x}]$.
\end{definition}

\subsection{Upper bounds and Laurent phenomenon}

\begin{definition}[upper bound $\U(\Sigma)$]
For a BFZ-seed $\Sigma$ its \emph{upper bound} is the $\Z\P$-subalgebra of $\F$ given by
\begin{equation*}
\U(\Sigma) = \Z\P\big[\mathbf{x}^{\pm 1}\big] \cap \Z\P\big[\mathbf{x_1}^{\pm 1}\big] \cap \dots \cap \Z\P\big[\mathbf{x_n}^{\pm 1}\big].
\end{equation*}
\end{definition}

The next theorem is a manifestation of the \emph{Laurent phenomenon} in terms of upper bounds.
\begin{theorem}[\protect{\cite[Theorem~1.5]{BFZ}}] \label{bfz-theorem}
If two seeds $\Sigma$ and $\Sigma'$ are related by a seed mutation and both are coprime,
then the corresponding upper bounds coincide: $\U(\Sigma) = \U(\Sigma')$.
\end{theorem}

\subsection{Relations between BFZ with \cite{GU} and this paper}

Given an exchange collection $V=(v_1,\dots,v_n)$ one can associate a skew-symmetric $n \times n$
matrix $B(V) = (b_{ij})$:
\begin{equation*} 
b_{i,j} = \omega(v_i, v_j).
\end{equation*}

\begin{lemma}
For any $V$ and $1\leq k\leq n$ we have $B(\mu_k V) = \mu_k B(V)$.
\end{lemma}

\begin{proof}
Indeed, let $B(V) = (b_{ij})$ and $B(\mu_k V) = (b'_{ij})$.
Then $b'_{ij} = \omega( v_i + \max(0,\omega(v_k,v_i)) v_k, v_j + \max(0,\omega(v_k,v_j)) v_k )
= \omega(v_i, v_j) + \max(0,\omega(v_k,v_i)) \omega(v_k,v_j) + \max(0,\omega(v_k,v_j)) \omega(v_i,v_k)
= b_{ij} +\max(0,-b_{ik}) b_{kj} + \max(0,b_{kj}) b_{ik} = b_{ij} + a \cdot b_{ik} b_{kj}$,
where $a=\frac{\operatorname{sgn}(b_{ik})+\operatorname{sgn}(b_{kj})}{2}$, i.e.\ $1$ if both $b_{ik}$ and $b_{kj}$ are positive,
$-1$~if they are both negative, and~$0$ otherwise.
It is easy to check that this coincides with Def\/inition~\ref{mub}.
\end{proof}

\begin{remark}
We note that in case $\rk L = 2$ the matrix $B(V)$ is a very special skew-symmetric matrix:
it is non-zero only if the collection~$V$ has at least two non-collinear vectors,
and in this case its rank equals two.
\end{remark}

\begin{remark}
For an exchange collection $V \in L^n$ the sublattice $L_V$ in $L$ denotes the sublattice generated by~$v_i$.
It can be seen that~$L_V$ is preserved under mutations of~$V$, and actually can be reconstructed from~$B(V)$
if~$\omega$ is non-degenerate on~$L_V$.
\end{remark}

\begin{remark} \label{compare-bfz}
Roughly the setup of BFZ corresponds to a special class of $C$-seeds
with $v_1,\dots,v_n$ being a base of the lattice~$L$ with all multiplicities
equal to~$1$. Thus the proof of Theorem~\ref{ublemma}
mostly reduces to Theorem~\ref{bfz-theorem} and its proof,
with extra care of keeping track of all the multiplicities
and exploiting nice functorial properties with respect to
the maps of the lattices and subcollections.

Lemmas~\ref{l41a} and~\ref{l42} are analogues
 and almost immediate consequences of Lemmas~4.1 and~4.2 in \cite{BFZ}.
Propositions~\ref{lhz2} and~\ref{l43c1} are analogues of Case $1$ in the proof of Proposition~4.3 in~\cite{BFZ}.
Proposition~\ref{lhz3} is analogue of the Case~2 in the proof of Proposition~4.3 of~\cite{BFZ}.
Lemma~\ref{l46} is similar to Lemma~4.6 in~\cite{BFZ}.
\end{remark}

\section{Def\/initions from \cite{GU}}\label{appendix2}

\begin{definition}[$U$-seed]
A $U$-seed is a quadruple $(W,V,F,X)$
where $V \in L_1^n$ is an exchange collection,
$F$ is a fan in $M$,
$X$ is a toric surface associated with the fan~$F$
and $W$ is a rational function on $X$. In addition, given a $U$-seed we can def\/ine a curve~$C$ by the equation
 \[C-\Sigma_{t}n_tD_t=(W),\]
where $\Sigma_{t} n_t D_t$ is the part corresponding to toric divisors.
\end{definition}

\begin{definition}[property $(U)$] \label{pu}
We say a $U$-seed satisf\/ies property $(U)$ if the following conditions hold:
\begin{enumerate}\itemsep=0pt
 \item[1)] $C$ is an ef\/fective divisor, i.e.\ $W$ is a Laurent polynomial;
 \item[2)] $C = A+B$, where $A$ is an irreducible non-rational curve and $B$ is supported on rational curves;
 \item[3)] the intersection of $C$ with toric divisors has canonical coordinates~$-1$;
 \item[4)] if $t \in V$, then the intersection index $(C \cdot D_t) \geq n_t$;
 \item[5)] for a toric divisor $D_t$ the intersection index $(A \cdot D_t)$ equals the number of $i$ such that $v_i = t$.
\end{enumerate}
\end{definition}

In \cite{GU} the Laurent phenomenon is established in the following terms

\begin{theorem}[$U$-lemma] \label{ulemma}
If two $U$-seeds $\Sigma$ and $\Sigma'$ are related by a mutation
then
$\Sigma$ satisfies property $(U)$
 $\iff$
$\Sigma'$ satisfies property $(U)$.
\end{theorem}

\subsection*{Acknowledgements}

We want to thank Bernhard Keller for interesting discussions
and in particular for suggesting the idea of considering upper bounds in the context of~\cite{GU}.
We thank Denis Auroux, Arkady Berenstein, Alexander Ef\/imov, Alexander Goncharov
and especially Alexandr Usnich
for ge\-ne\-rous\-ly sharing their insights on mutations or mirror symmetry.
We would like to thank Martin Guest for reading the draft of this paper and his valuable comments.
We also want to thank the anonymous referees for their helpful comments and suggestions.
First author was supported by a~Japanese Government (Monbukagakusho:MEXT) Scholarship. Also this work was supported by
World Premier International Research Center Initiative (WPI Initiative), MEXT, Japan,
Grant-in-Aid for Scientif\/ic Research (10554503) from Japan Society for Promotion of Science,
Grant of Leading Scientif\/ic Schools (N.Sh.\ 4713.2010.1), grants ERC GEMIS and FWF P 24572-N25.
The authors thank the IPMU for its hospitality and support during their visits in May
and November 2012 which helped to f\/inish this work.

\pdfbookmark[1]{References}{ref}
\LastPageEnding

\end{document}